\documentclass[leqno,12pt]{amsart}
\usepackage{amsmath,amsthm}
\usepackage{amscd}
\usepackage[psamsfonts]{amssymb}
\usepackage[mathscr,psamsfonts]{eucal}
\DeclareMathAlphabet{\eurm}{U}{eur}{m}{n}

\newcommand{\E}[1]{{\eurm{#1}}}

\newcommand{\f}[1]{{\boldsymbol{#1}}}
\DeclareMathOperator{\id}{{id}}

\DeclareMathOperator{\Alt}{{Alt}}

\DeclareMathOperator{\byd}{\,{\raisebox{.1ex}{$\eurm :$}{\eurm =}}\,}
\newcommand{\sig}{\sigma}

\newcommand{\del}{\delta}
\newcommand{\lam}{\lambda}
\newcommand{\ome}{\omega}
\newcommand{\Lam}{\Lambda}
\newcommand{\Gam}{\Gamma}

\newcommand{\bEq}{\begin{eqnarray}}
\newcommand{\eEq}{\end{eqnarray}}
\newcommand{\beq}{\begin{eqnarray*}}
\newcommand{\eeq}{\end{eqnarray*}}
\newcommand{\QED}{{\,\text{\footnotesize QED}}}

\newcommand{\dt}[1]{{\dot{#1}}}
\newcommand{\uten}[1]{\underset{#1}{\otimes}}
\newcommand{\ucar}[1]{\underset{#1}{\times}}
\newcommand{\sep}[1]{\qquad\text{#1}\qquad}

\newcommand{\col}[3]{_{#1}{}^{#2}{}_{#3}}
\newcommand{\der}{\partial}
\newcommand{\nab}{\nabla}

\newcommand{\lang}{\langle}
\newcommand{\rang}{\rangle}

\newcommand{\com}{\circ}
\newcommand{\car}{\times}
\newcommand{\ten}{\otimes}
\newcommand{\wed}{\wedge}
\newcommand{\Rn}{{I\!\!R}}

\newcommand{\END}{{\,\text{\qedsymbol}}}

\begin{document}

\newcounter{theorem}

\newtheorem{definition}[theorem]{Definition}
\newtheorem{lemma}[theorem]{Lemma}
\newtheorem{proposition}[theorem]{Proposition}
\newtheorem{theorem}[theorem]{Theorem}
\newtheorem{corollary}[theorem]{Corollary}
\newtheorem{remark}[theorem]{Remark}
\newtheorem{example}[theorem]{Example}
\newtheorem{Note}[theorem]{Note}
\newcounter{assump}
\newtheorem{Assumption}{\indent Assumption}[assump]
\renewcommand{\thetheorem}{\thesection.\arabic{theorem}}

\newcommand{\bCr}{\begin{corollary}}
\newcommand{\eCr}{\end{corollary}}
\newcommand{\bDf}{\begin{definition}\em}
\newcommand{\eDf}{\end{definition}}
\newcommand{\bLm}{\begin{lemma}}
\newcommand{\eLm}{\end{lemma}}
\newcommand{\bPr}{\begin{proposition}}
\newcommand{\ePr}{\end{proposition}}
\newcommand{\bRm}{\begin{remark}\em}
\newcommand{\eRm}{\end{remark}}
\newcommand{\bEx}{\begin{example}\em}
\newcommand{\eEx}{\end{example}}
\newcommand{\bTh}{\begin{theorem}}
\newcommand{\eTh}{\end{theorem}}
\newcommand{\bNt}{\begin{Note}\em}
\newcommand{\eNt}{\end{Note}}
\newcommand{\bPf}{\begin{proof}[\noindent\indent{\sc Proof}]}
\newcommand{\ePf}{\renewcommand{\qedsymbol}{}\end{proof}}
\title[On the curvature of tensor product connections]
{On the curvature of tensor product connections\\
        and covariant differentials}

\author[Josef Jany\v ska]
{Josef Jany\v ska}

\address
{
\newline
Department of Mathematics, Masaryk University
\newline
Jan\'a\v ckovo n\'am. 2a, 662 95 Brno,
Czech Republic
\newline
E-mail: {\tt janyska@math.muni.cz}
}

\thanks
{This paper is in final form and no version of it will be submitted
for publication elsewhere.
\\
\indent
This paper has been supported
by the Grant agency of the Czech Republic under the project
number GA 201/02/0225.
}



\maketitle
\begin{abstract}
We give coordinate formula and geometric description
of the curvature of the tensor product connection of
linear connections on vector bundles with the same base manifold. We
define the covariant differential of geometric fields of certain types
with respect to a pair of a linear connection on a vector bundle
and a linear symmetric connection on the base manifold. We prove
the generalized Bianchi identity for linear connections and we prove
that the antisymmetrization of the second order covariant differential
is expressed via the curvature tensors of both
connections.
\end{abstract}

\section*{Introduction}

In the theory of linear symmetric (classical) connections on
a manifold there are many very well known identities of the curvature
tensor (see for instance \cite{KobNom63,Sch54}).
Some of these identities can be generalized for any linear
connection on a vector bundle.

In this paper we give the coordinate
formula for the curvature of the tensor product connection
$K\ten K'$ of two linear connections $K$ or $K'$ on vector bundles
$\f E\to \f M$ or $\f E'\to \f M$, respectively, and we give also
the geometric description of this curvature. We prove that
the curvature of $K\ten K'$ is determined by the curvatures
of $K$ and $K'$.

The above results are used in the case if one of linear connections
is a classical (linear and symmetric) connection
on the base manifold. We introduce the
covariant differential of sections of tensor products
(over the base manifold) of a vector bundle, its dual vector bundle,
the tangent and the cotangent bundles of the base manifold. We
prove that such (first order) covariant differential of
the curvature tensor of a linear connection satisfies the generalized
Bianchi identity and that the antisymmetrization of the second order
covariant differential is expressed through the curvatures
of linear and classical connections.

All manifolds and maps are supposed to be smooth.

\section{Linear connections on vector bundles}
\label{Linear connections on vector bundles}
\setcounter{equation}{0}
\setcounter{theorem}{0}

Let $p:\f E\to \f M$ be a vector bundle.
Local linear fiber coordinate
charts on $\f E$ will be denoted by $(x^\lam, y^i)$.
The corresponding base of local sections of $\f E$ or $\f E^*$
will be denoted by $\E b_i$ or $\E b^i$, respectively.

\smallskip

\bDf\label{Df1.1}
We define a {\em linear connection\/} on $\f E$ to be a linear
splitting
\beq
K : \f E \to J^1\f E\,.\quad\END
\eeq
\eDf

\bPr\label{Pr1.2}
Considering the contact morphism $J^1\f E\to
T^*\f M \ten T\f E$ over the identity of $T\f M$, a
linear connection can be regarded as a $T\f E$-valued
1-form
\beq
K : \f E \to T^*\f M \ten T\f E
\eeq
projecting on the identity of $T\f M$.

The coordinate expression of a linear connection
$K$
is of the form
\beq
K = d^\lam \ten \big(\der_\lam
+ K\col ji\lam \, y^j \, \der_i \big) \,,
\sep{with}  K\col ji \lam \in C^\infty(\f M,\Rn)\,. \quad\END
\eeq
\ePr

\bDf\label{Df1.3}
The {\em covariant differential} of a section $\Phi:\f M\to \f E$
with respect to $K$ is defined to be
\beq
\nab^K\Phi = j^1\Phi - K\com\Phi: \f M \to \f E\uten{\f M}T^*\f M\,.
\quad\END
\eeq
\eDf

\bRm\label{Rm1.4}
From the affine structure of $\pi^1_0:J^1\f E\to \f E$
we obtain that the difference
$j^1\Phi - K\com\Phi$ lies in the associated vector bundle
$V\f E\ten T^*\f M$. From $V\f E=\f E\ucar{\f M}\f E$
we get the above Definition \ref{Df1.3}.
\END
\eRm

Let $\Phi=\phi^i \, \E b_i$, then we have the coordinate expression
\beq
\nab^K\Phi = (\der_\lam \phi^i - K\col ji\lam \phi^j)\,
        \E b_i\ten d^\lam\,.
\eeq

\bDf\label{Df1.5}
The {\em curvature} of a linear connection
${K}$
on
$\f E$
turns out to be the vertical valued 2--form
\beq
R[{K}] = - [{K},{K}]:\f E\to V\f E\ten \Lam^2 T^*\f M\,,
\eeq
where $[,]$ is the Froelicher-Nijenhuis bracket.
\END
\eDf

The coordinate
expression is
\begin{align*}
R[{K}] & = R[K]\col ji{\lam\mu}\, y^j\,  \der_i \ten d^\lam \wed d^\mu
\\
& = -2(\der_\lam K \col ji\mu + K\col jp\lam K\col pi\mu)\, y^j\,
        \der_i \ten d^\lam \wed d^\mu \,,
\end{align*}
i.e. the coefficients of the curvature are
\beq
R[{K}]\col ji{\lam\mu} = \der_\mu K \col ji\lam -
        \der_\lam K \col ji\mu  +
        K\col jp\mu K\col pi\lam - K\col jp\lam K\col pi\mu
        \,.
\eeq

If we consider the identification $V\f E = \f E \ucar{\f M}\f E$
and linearity of $R[K]$, the curvature $R[K]$ can be considered as
a tensor field (the curvature tensor field)
$
R[{K}] :\f M\to\f E^* \ten \f E \ten \Lam^2 T^*\f M\,.
$

\smallskip

\bTh\label{Th1.6}
We have the generalized Bianchi identity
\beq
 [K, R[{K}]]=0\,.
\eeq
\eTh

\bPf
It follows from the graded Jacobi identity for the Froelicher-Nijenhuis
bracket. \QED
\ePf

We have, \cite{ManMod89},

\bPr\label{Pr1.7}
Let $K$ be a linear connection on $\f E$. Then, there is a unique
linear connection $K^*: \f E^*\to J^1\f E^*$ on the dual
vector bundle $\f E^*\to \f M$ such that the following diagram
commutes
\begin{equation*}\begin{CD}
\f E\ucar{\f M}\f E^*
        @>{\lang,\rang}>>
        \f M\car \Rn
\\
        @V{K\car K^*}VV  @V{0\car \id_{\Rn}}VV
\\
        J^1\f E\ucar{\f M}J^1\f E^*
        @>{J^1\lang,\rang}>>
        T^*\f M \car \Rn
\end{CD}\end{equation*}

Its coordinate expression is
\beq
K^* = d^\lam \ten \big(\der_\lam
- K\col ij\lam \, y_j \, \der^i \big) \,,
\sep{with}  K\col ij \lam \in C^\infty(\f M,\Rn)\,,
\eeq
where $(x^\lam,y_i)$ are the induced linear fiber
coordinates on $\f E^*$
and $\der^i=\der/\der y_i$. \END
\ePr

\bDf\label{Df1.8}
The connection $K^*$ is said to be the
{\em dual connection of $K$}.
\END
\eDf

\bPr\label{Pr1.8}
We have
$
R[K^*]:\f M\to \f E\uten{\f M}\f E^*\uten{\f M}\Lam^2T^*\f M
$
and
\beq
R[K^*]^i{}_{j\lam\mu}=-R[K]\col ji{\lam\mu}\,.\quad\END
\eeq
\ePr

\section{Tensor product linear connections}
\label{Tensor product linear connections}
\setcounter{equation}{0}
\setcounter{theorem}{0}

Let $p':\f E'\to \f M$ be another vector bundle.
Local linear fiber coordinate
charts on $\f E'$ will be denoted by $(x^\lam, z^a)$.
The corresponding base of local sections of $\f E'$
or $\f E'{}^*$
will be denoted by $\E b'_a$ or $\E b'{}^a$, respectively.

\smallskip

Consider a linear connection $K'$ on $\f E'$
with coordinate expression
\beq
K' = d^\lam \ten \big(\der_\lam
+ K'\col ba\lam \, z^b \, \der_a \big) \,,
\sep{with}  K'\col ba \lam \in C^\infty(\f M,\Rn)\,.
\eeq

\smallskip
Let us consider the tensor product $\f E\uten{\f M}
\f E'\to \f M$
with the induced fiber linear coordinate chart $(x^\lam, w^{ia})$.
We have, \cite{ManMod89},

\bPr\label{Pr2.1}
Let $K$ be a linear connection on $\f E$ and $K'$ be a
linear connection
on $\f E'{}$. Then, there is a unique linear connection
$K\ten K':\f E\uten{\f M} \f E'{}\to
J^1(\f E\uten{\f M} \f E')$ such that the following diagram commutes
\begin{equation*}\begin{CD}
\f E\ucar{\f M}\f E'
        @>{\ten}>>
        \f E\uten{\f M}\f E'
\\
        @V{K\car K'}VV  @V{K\ten K'}VV
\\
        J^1\f E\ucar{\f M}J^1\f E'
        @>{J^1\ten}>>
        J^1(\f E\uten{\f M}\f E')
\end{CD}\end{equation*}

Its coordinate expression is
$$
K\ten K' = d^\lam\ten\big(\der_\lam + (K\col ji\lam w^{ja}
        + K'\col ba\lam w^{ib})\der_{ia}\big)\,.\quad\END
$$
\ePr

\bDf\label{Df2.2}
The connection $K\ten K'$ is said to be the {\em tensor product
connection} of $K$ and $K'$. \END
\eDf

\bRm
We remark that this concept was introduced in another way in
\cite{KolMicSlo93}, p. 381.
\END\eRm

The tensor product connection is linear, so we can define its
tensor product connection with another linear connection and we have
by iteration

\bPr\label{Pr2.3}
A linear connection $K$ on $\f E$ and a linear connection $K'$
on $\f E'$ induce the linear tensor product connection
$K^p_q \ten K'{}^r_s\byd
\ten^p K \ten\ten^q K^*\ten\ten^r K'\ten\ten^s K'{}^*$
on $\ten^p\f E\uten{\f M} \ten^q\f E^*\uten{\f M}
\ten^r\f E'\uten{\f M} \ten^s\f E'{}^*$
with coordinate expression
\begin{align*}
K^p_q\ten K'{}^r_s & =d^\lam\ten\bigg(\der_\lam +
     (K\col {k}{i_1}\lam \,
     w^{ki_2\dots i_p{a}_1\dots{a}_r}_{j_1\dots j_q{b}_1\dots{b}_s}
     +\dots + K\col {k}{i_p}\lam \,
     w^{i_1\dots i_{p-1}k{a}_1\dots{a}_r}_{j_1\dots
                j_q{b}_1\dots{b}_s}
\\
 & \quad -  K\col {j_1}{k}\lam \,
     w^{i_1\dots i_p{a}_1\dots{a}_r}_{kj_2\dots j_q{b}_1\dots{b}_s}
     - \dots - K\col {j_q}{k}\lam \,
     w^{i_1\dots i_{p}{a}_1\dots{a}_r}_{j_1\dots j_{q-1}k{b}_1
                \dots{b}_s}
        \nonumber
\\
  & \quad + {K'}\col {{c}}{{a}_1}\lam \,
     w^{i_1\dots i_p{c}{a}_2\dots{a}_r}_{j_1
                \dots j_q{b}_1\dots{b}_s}
     + \dots + {K'}\col {{c}}{{a}_r}\lam \,
     w^{i_1\dots i_{p}{a}_1\dots{a}_{r-1}{c}}_{j_1
                \dots j_q{b}_1\dots{b}_s}
        \nonumber
\\
 & \quad -  {K'}\col {{b}_1}{{c}}\lam \,
     w^{i_1\dots i_p{a}_1\dots{a}_r}_{j_1\dots j_q{c}{b}_2
                \dots{b}_s}
     - \dots - {K'}\col {{b}_s}{{c}}\lam \,
     w^{i_1\dots i_{p}{a}_1\dots{a}_r}_{j_1\dots j_{q}{b}_1
                \dots{b}_{s-1}{c}}
     )\, \der_{i_1\dots i_p{a}_1\dots {a}_r}^{j_1\dots j_q{b}_1
                \dots{b}_s}\bigg)
        \nonumber
\end{align*}
where
$(x^\lam, w^{i_1\dots i_pa_1\dots a_r}_{j_1\dots j_q b_1
        \dots b_s})$
are the induced linear fiber coordinates
on $\ten^p\f E\uten{\f M} \ten^q\f E^*\uten{\f M}
\ten^r\f E'\uten{\f M} \ten^s\f E'{}^*$.
\END
\ePr

\smallskip

The curvature of the linear tensor product connection
${K\ten K'}$
on
$\f E\uten{\f M}\f E'$
turns out to be the vertical valued 2--form
\beq
R[{K\ten K'}] = - [{K\ten K'},{K\ten K'}]:
\f E\uten{\f M}\f E'\to V(\f E\uten{\f M}\f E')
        \ten \Lam^2 T^*\f M\,.
\eeq

\bTh\label{Th2.4}
The coordinate
expression of $R[K\ten K']$ is
\begin{align*}
R[{K\ten K'}] & = R[K\ten K']\col {jb}{ia}{\lam\mu}\, w^{jb}\,
        \der_{ia} \ten d^\lam \wed d^\mu
\\
& = \big(R[K] \col ji{\lam\mu} \, w^{ja} +
        R[K']\col ba{\lam\mu} \, w^{ib}\big)
        \der_{ia} \ten d^\lam \wed d^\mu \,,
\end{align*}
i.e. the coefficients of the curvature $R[K\ten K']$ are
\beq
R[K\ten {K'}]\col {jb}{ia}{\lam\mu}= R[K]\col ji{\lam\mu}\del^a_b
        + R[K']\col ba{\lam\mu}\del^i_j\,.
\eeq
\eTh

\bPf
This can be proved in coordinates.
\QED\ePf

Theorem \ref{Th2.4} implies that the curvature $R[K\ten K']$
is determined by the curvatures $R[K]$ and $R[K']$.
Now, we would like to find the geometric description of the
curvature $R[K\ten K']$. First we note that
the curvatures of the above linear connections
can be considered as
bilinear morphisms, over $\f M$,
\begin{align*}
R[K] & : \f E\ucar{\f M} \f E^* \to \Lam^2T^*\f M\,,
\\
R[K'] & : \f E'\ucar{\f M} \f E'{}^* \to \Lam^2T^*\f M\,,
\\
R[K\ten K'] & : (\f E\uten{\f M}\f E')\ucar{\f M}
        (\f E\uten{\f M}\f E')^* \to \Lam^2T^*\f M\,.
\end{align*}

Then we have

\bTh\label{Th2.5}
The curvature $R[K\ten K']$ is a unique bilinear morphism
such that the following diagram commutes
\begin{equation*}\begin{CD}
\f E\ucar{\f M}\f E'\ucar{\f M}\f E^*\ucar{\f M}\f E'{}^*
        @>{\lang,\rang' R[K]+\lang,\rang R[K']}>>
        \Lam^2 T^*\f M
\\
        @V{(\ten,\ten)}VV  @V{\id_{\Lam^2 T^*\f M}}VV
\\
        (\f E\uten{\f M}\f E')\ucar{\f M}(\f E^*\uten{\f M}\f E'{}^*)
        @>{R[K\ten K']}>>
        \Lam^2 T^*\f M
\end{CD}\end{equation*}
where $\lang,\rang$ or $\lang,\rang'$ are the evaluation morphisms
on $\f E$ or $\f E'$, respectively.
\eTh

\bPf
Let us assume a bilinear morphism
$R: (\f E\uten{\f M}\f E')\ucar{\f M}
        (\f E\uten{\f M}\f E')^* \to \Lam^2T^*\f M$
and let us put $e=(e^i)\in \f E_x$, $e^*=(e_i)\in \f E^*_x$,
$e'=(e'{}^a)\in \f E'_x{}$ and $e'{}^*=(e'_a)\in \f E'_x{}^*$.
Then
\begin{align*}
\lang e',e'{}^*\rang\, R[K] (e, e^*)
        & = e'{}^a e'_a R[K]\col ji{\lam\mu} e^j e_i
        \, d^\lam\wed d^\mu\,,
\\
\lang e,e^*\rang\, R[K'] (e', e'{}^*)
        & = e^i e_i R[K']\col ba{\lam\mu} e'{}^b e'_a
        \, d^\lam\wed d^\mu\,,
\\
R (e\ten e', e^*\ten e'{}^*)
        & = R\col {jb}{ia}{\lam\mu} e^j e'{}^b e_ie'_a
        \, d^\lam\wed d^\mu\,
\end{align*}
and it is easy to see that
$R (e\ten e', e^*\ten e'{}^*)= \lang e',e'{}^*\rang\, R[K] (e, e^*)
+ \lang e,e^*\rang\, R[K'] (e', e'{}^*)$
if and only if
\beq
R\col {jb}{ia}{\lam\mu} = R[K]\col ji{\lam\mu} \del^a_b
        + R[K']\col ba{\lam\mu} \del^i_j \,.
\eeq
Now, Theorem \ref{Th2.5} follows from Theorem \ref{Th2.4}.
\QED
\ePf

\bPr\label{Pr2.6}
The curvature
$
R[K^p_q\ten {{K'}}^r_s]\byd - [K^p_q\ten {{K'}}^r_s , K^p_q
                \ten {{K'}}^r_s]
$
is determined by the curvatures $R[K]$ and $R[{{K'}}]$.
We have the coordinate expression
\begin{align*}
R[K^p_q\ten {{K'}}^r_s] & = \bigg(R[K]\col {k}{i_1}{\lam\mu} \,
     w^{ki_2\dots i_p{a}_1\dots{a}_r}_{j_1\dots j_q{b}_1\dots{b}_s}
     +\dots + R[K]\col {k}{i_p}{\lam\mu} \,
     w^{i_1\dots i_{p-1}k{a}_1\dots{a}_r}_{j_1\dots j_q{b}_1
        \dots{b}_s}
\\
 & \quad -  R[K]\col {j_1}{k}{\lam\mu} \,
     w^{i_1\dots i_p{a}_1\dots{a}_r}_{kj_2\dots j_q{b}_1\dots{b}_s}
     - \dots - R[K]\col {j_q}{k}{\lam\mu} \,
     w^{i_1\dots i_{p}{a}_1\dots{a}_r}_{j_1\dots j_{q-1}k{b}_1
        \dots{b}_s}
        \nonumber
\\
  & \quad + R[{K'}]\col {{c}}{{a}_1}{\lam\mu} \,
     w^{i_1\dots i_p{c}{a}_2\dots{a}_r}_{j_1\dots j_q{b}_1
        \dots{b}_s}
     + \dots + R[{K'}]\col {{c}}{{a}_r}{\lam\mu} \,
     w^{i_1\dots i_{p}{a}_1\dots{a}_{r-1}{c}}_{j_1\dots j_q{b}_1
        \dots{b}_s}
        \nonumber
\\
 & \quad -  R[{K'}]\col {{b}_1}{{c}}{\lam\mu} \,
     w^{i_1\dots i_p{a}_1\dots{a}_r}_{j_1\dots j_q{c}{b}_2
        \dots{b}_s}
     - \dots - R[{K'}]\col {{b}_s}{{c}}{\lam\mu} \,
     w^{i_1\dots i_{p}{a}_1\dots{a}_r}_{j_1\dots j_{q}{b}_1
        \dots{b}_{s-1}{c}}
     \bigg)\,
\\
 &\qquad \E b_{i_1\dots i_p}\ten\E b^{j_1 \dots j_q}\ten
        \E b_{{a}_1\dots{a}_r} \ten \E b^{{b}_1\dots{b}_s}
        \ten d^{\lam}\wed d^{\mu}
        \,,
\end{align*}
where we have put $\E b_{i_1\dots i_p}=
\E b_{i_1}\ten \dots \ten \E b_{i_p}$,
$\E b^{j_1 \dots j_q} = \E b^{j_1}\ten
        \dots \ten \E b^{j_q}$,
$\E b_{{a}_1\dots{a}_r} =
\E b_{{a}_1}\ten\dots\ten \E b_{{a}_r}$,
$ \E b^{{b}_1\dots{b}_s} =
        \E b ^{{b}_1}\ten\dots\ten \E b^{{b}_s}$.
\ePr

\bPf
This follows from the definition of the curvature,
Proposition \ref{Pr1.8} and the iteration
of Theorem \ref{Th2.4}.
\END
\ePf

\section{Classical connections}
\label{Classical connections}
\setcounter{equation}{0}
\setcounter{theorem}{0}

Let $\f M$ be an $m$-dimensional manifold.
Local coordinate charts on $\f M$ will be denoted by
$(x^\lam)$, $\lam=1,\dots,m$, the induced coordinate charts on
$T\f M$ or $T^*\f M$
will be denoted by
$(x^\lam, \dot x^\lam)$ or $(x^\lam, \dot x_\lam)$ and the induced
local bases of sections of
$T\f M$ or
$T^*\f M$
are denoted by
$(\der_\lam)$ or
$(d^\lam)$, respectively.

\smallskip

A {\em classical connection\/} on $\f M$ is defined to be
a linear symmetric connection on $p_{\f M}:T\f M\to \f M$
with coordinate expression
\beq
\Gam = d^\lam \ten \big(\der_\lam
+ \Gam\col\nu\mu\lam \, \dot x^\nu \, \dt\der_\mu \big) \,,
\quad  \Gam\col \mu \lam \nu\in C^\infty(\f M,\Rn),\quad
\Gam\col \mu \lam \nu = \Gam\col \nu \lam \mu \,.
\eeq

\bRm\label{Rm3.3}
Let us recall the 1st and the 2nd Bianchi identities of classical
connections expressed in coordinates by
\begin{align*}
R[\Gam]\col {\nu}\rho{\lam\mu} + R[\Gam]\col {\lam}\rho{\mu\nu}
        + R[\Gam]\col {\mu}\rho{\nu\lam} & = 0\,,
\\
R[\Gam]\col {\nu}\rho{\lam\mu;\sig}+
R[\Gam]\col {\nu}\rho{\mu\sig;\lam}+
R[\Gam]\col {\nu}\rho{\sig\lam;\mu} & = 0\,,
\end{align*}
respectively, where $;$ denotes the covariant differential
with respect to $\Gam$.
\END
\eRm

\smallskip
Let us denote by
$
\f E^{p,r}_{q,s}\byd \ten^p \f E
\uten{\f M}\ten^q \f E^*\uten{\f M}\ten^rT\f M
\uten{\f M} \ten^s T^*\f M
$. Then, as a direct consequence of Proposition \ref{Pr2.3},
we have

\bPr\label{Pr3.4}
A classical connection $\Gam$ on $\f M$ and a linear
connection ${K}$ on $\f E$ induce the linear tensor product connection
$K^p_q\ten {\Gam}^r_s\byd \ten^p K\ten\ten^q K^*\ten
\ten^r{\Gam} \ten \ten^s{\Gam}^*$ on
$\f E^{p,r}_{q,s}$
\beq
K^p_q\ten {\Gam}^r_s: \f E^{p,r}_{q,s} \to T^*\f M\uten{\f M}T
        \f E^{p,r}_{q,s}
\eeq
with coordinate expression
\begin{align}\nonumber
K^p_q\ten {\Gam}^r_s & =d^\nu\ten\bigg(\der_\nu +
     (K\col {k}{i_1}\nu \,
     y^{ki_2\dots i_p\lam_1\dots\lam_r}_{j_1\dots j_q\mu_1\dots\mu_s}
     +\dots + K\col {k}{i_p}\nu \,
     y^{i_1\dots i_{p-1}k\lam_1\dots\lam_r}_{j_1\dots
                j_q\mu_1\dots\mu_s}
\\
 & \quad -  K\col {j_1}{k}\nu \,
     y^{i_1\dots i_p\lam_1\dots\lam_r}_{kj_2\dots j_q\mu_1\dots\mu_s}
     - \dots - K\col {j_q}{k}\nu \,
     y^{i_1\dots i_{p}\lam_1\dots\lam_r}_{j_1\dots j_{q-1}k\mu_1
                \dots\mu_s}
        \nonumber
\\
  & \quad + \Gam\col {\rho}{\lam_1}\nu \,
     y^{i_1\dots i_p\rho\lam_2\dots\lam_r}_{j_1
                \dots j_q\mu_1\dots\mu_s}
     + \dots + \Gam\col {\rho}{\lam_r}\nu \,
     y^{i_1\dots i_{p}\lam_1\dots\lam_{r-1}\rho}_{j_1
                \dots j_q\mu_1\dots\mu_s}
        \nonumber
\\
 & \quad -  \Gam\col {\mu_1}{\rho}\nu \,
     y^{i_1\dots i_p\lam_1\dots\lam_r}_{j_1\dots j_q\rho\mu_2
                \dots\mu_s}
     - \dots - \Gam\col {\mu_s}{\rho}\nu \,
     y^{i_1\dots i_{p}\lam_1\dots\lam_r}_{j_1\dots j_{q}\mu_1
                \dots\mu_{s-1}\rho}
     )\, \der_{i_1\dots i_p\lam_1\dots \lam_r}^{j_1\dots j_q\mu_1
                \dots\mu_s}\bigg)
        \nonumber
\end{align}
where
$(x^\lam, y^{i_1\dots i_p\lam_1\dots \lam_r}_{j_1\dots j_q\mu_1
        \dots\mu_s})$
are the induced linear fiber coordinates
on $\f E^{p,r}_{q,s}$.
\END
\ePr

As a direct consequence of Proposition \ref{Pr2.6} we have

\bPr\label{Pr3.5}
The curvature
$
R[K^p_q\ten {\Gam}^r_s]
$
is determined by the curvatures $R[K]$ and $R[{\Gam}]$.
We have the coordinate expression
\begin{align*}
R[K^p_q\ten {\Gam}^r_s] & = \bigg(R[K]\col {k}{i_1}{\nu_1\nu_2} \,
     y^{ki_2\dots i_p\lam_1\dots\lam_r}_{j_1\dots j_q\mu_1\dots\mu_s}
     +\dots + R[K]\col {k}{i_p}{\nu_1\nu_2} \,
     y^{i_1\dots i_{p-1}k\lam_1\dots\lam_r}_{j_1\dots j_q\mu_1
        \dots\mu_s}
\\
 & \quad -  R[K]\col {j_1}{k}{\nu_1\nu_2} \,
     y^{i_1\dots i_p\lam_1\dots\lam_r}_{kj_2\dots j_q\mu_1\dots\mu_s}
     - \dots - R[K]\col {j_q}{k}{\nu_1\nu_2} \,
     y^{i_1\dots i_{p}\lam_1\dots\lam_r}_{j_1\dots j_{q-1}k\mu_1
        \dots\mu_s}
        \nonumber
\\
  & \quad + R[\Gam]\col {\rho}{\lam_1}{\nu_1\nu_2} \,
     y^{i_1\dots i_p\rho\lam_2\dots\lam_r}_{j_1\dots j_q\mu_1
        \dots\mu_s}
     + \dots + R[\Gam]\col {\rho}{\lam_r}{\nu_1\nu_2} \,
     y^{i_1\dots i_{p}\lam_1\dots\lam_{r-1}\rho}_{j_1\dots j_q\mu_1
        \dots\mu_s}
        \nonumber
\\
 & \quad -  R[\Gam]\col {\mu_1}{\rho}{\nu_1\nu_2} \,
     y^{i_1\dots i_p\lam_1\dots\lam_r}_{j_1\dots j_q\rho\mu_2
        \dots\mu_s}
     - \dots - R[\Gam]\col {\mu_s}{\rho}{\nu_1\nu_2} \,
     y^{i_1\dots i_{p}\lam_1\dots\lam_r}_{j_1\dots j_{q}\mu_1
        \dots\mu_{s-1}\rho}
     \bigg)\,
\\
 &\qquad \E b_{i_1\dots i_p}\ten\E b^{j_1 \dots j_q}\ten
        \der_{\lam_1\dots\lam_r} \ten d^{\mu_1\dots\mu_s}
        \ten d^{\nu_1}\wed d^{\nu_2}
        \,,
\end{align*}
where we have put $\E b_{i_1\dots i_p}=
\E b_{i_1}\ten \dots \ten \E b_{i_p}$,
$\E b^{j_1 \dots j_q} = \E b^{j_1}\ten
        \dots \ten \E b^{j_q}$,
$\der_{\lam_1\dots\lam_r} =
\der_{\lam_1}\ten\dots\ten \der_{\lam_r}$,
$ d^{\mu_1\dots\mu_s} =
        d^{\mu_1}\ten\dots\ten d^{\mu_s}$.
\END
\ePr

\section{Covariant differentials}
\label{Covariant differentials}
\setcounter{equation}{0}
\setcounter{theorem}{0}

Let us note that the tensor product connection
$K^p_q\ten{\Gam}^r_s$ can be considered
as a linear splitting
\beq
K^p_q\ten{\Gam}^r_s:\f E^{p,r}_{q,s}\to J^1\f E^{p,r}_{q,s}\,.
\eeq

\bDf\label{Df4.1}
Let $\Phi\in C^\infty(\f E^{p,r}_{q,s})$.
We define the {\em covariant differential of $\Phi$ with respect to
a pair of connections $(K,\Gam)$} as a section of
$\f E^{p,r}_{q,s}\ten T^*\f M$ defined by
\beq
\nab^{(K,\Gam)} \Phi = j^1 \Phi - (K^p_q\ten{\Gam}^r_s )\com \Phi\,.
\quad\END
\eeq
\eDf

\bRm\label{Rm4.2}
The covariant differential $\nabla^{(K,\Gam)}\Phi$ is in fact
the standard covariant differential (see Definition \ref{Df1.3})
$\nabla^{K^p_q\ten \Gam^r_s}\Phi$.
\END
\eRm

\bPr\label{Pr4.3}
Let $\Phi\in C^\infty(\f E^{p,r}_{q,s})$,
        $\Phi=\phi^{i_1\dots i_p\lam_1\dots\lam_r}
        _{j_1 \dots j_q\mu_1\dots\mu_s}\,
        \E b_{i_1\dots i_p}\ten\E b^{j_1 \dots j_q}\ten
        \der_{\lam_1\dots\lam_r} \ten d^{\mu_1\dots\mu_s}$.
Then we have the coordinate expression
\begin{align*}
\nab^{(K,\Gam)} \Phi & =\nab^{(K,\Gam)}
        _{\nu}\phi^{i_1\dots i_p\lam_1\dots\lam_r}
        _{j_1 \dots j_q\mu_1\dots\mu_s}\,
        \E b_{i_1\dots i_p}\ten\E b^{j_1 \dots j_q}\ten
        \der_{\lam_1\dots\lam_r} \ten d^{\mu_1\dots\mu_s}\ten d^\nu
\\
& = \bigg(\der_{\nu}\phi^{i_1\dots i_p\lam_1\dots\lam_r}
        _{j_1 \dots j_q\mu_1\dots\mu_s}
     -K\col {k}{i_1}\nu \,
     \phi^{ki_2\dots i_p\lam_1\dots\lam_r}_{j_1\dots j_q\mu_1\dots\mu_s}
     -\dots - K\col {k}{i_p}\nu \,
     \phi^{i_1\dots i_{p-1}k\lam_1\dots\lam_r}_{j_1\dots
                j_q\mu_1\dots\mu_s}
\\
 & \quad +  K\col {j_1}{k}\nu \,
     \phi^{i_1\dots i_p\lam_1\dots\lam_r}_{kj_2\dots j_q\mu_1\dots\mu_s}
     + \dots + K\col {j_q}{k}\nu \,
     \phi^{i_1\dots i_{p}\lam_1\dots\lam_r}_{j_1\dots j_{q-1}k\mu_1
                \dots\mu_s}
        \nonumber
\\
  & \quad - \Gam\col {\rho}{\lam_1}\nu \,
     \phi^{i_1\dots i_p\rho\lam_2\dots\lam_r}_{j_1
                \dots j_q\mu_1\dots\mu_s}
     - \dots - \Gam\col {\rho}{\lam_r}\nu \,
     \phi^{i_1\dots i_{p}\lam_1\dots\lam_{r-1}\rho}_{j_1
                \dots j_q\mu_1\dots\mu_s}
        \nonumber
\\
 & \quad +  \Gam\col {\mu_1}{\rho}\nu \,
     \phi^{i_1\dots i_p\lam_1\dots\lam_r}_{j_1\dots j_q\rho\mu_2
                \dots\mu_s}
     + \dots + \Gam\col {\mu_s}{\rho}\nu \,
     \phi^{i_1\dots i_{p}\lam_1\dots\lam_r}_{j_1\dots j_{q}\mu_1
                \dots\mu_{s-1}\rho}
     \bigg)
\\
& \qquad \E b_{i_1\dots i_p}\ten\E b^{j_1 \dots j_q}\ten
        \der_{\lam_1\dots\lam_r} \ten d^{\mu_1\dots\mu_s}
        \ten d^{\nu}\,.
\end{align*}
\ePr

\bPf
The proof follows immediately from Definition \ref{Df4.1}
and the coordinate expression (see Proposition \ref{Pr3.4})
of the connection
$K^p_q\ten{\Gam}^r_s$. \QED
\ePf

In what follows we set $\nab = \nab^{(K,\Gam)}$ and
$\phi^{i_1\dots i_p\lam_1\dots\lam_r}
        _{j_1 \dots j_q\mu_1\dots\mu_s;\nu}=
\nab_{\nu}\phi^{i_1\dots i_p\lam_1\dots\lam_r}
        _{j_1 \dots j_q\mu_1\dots\mu_s}$.

\bRm
If $p=q=0$ the field $\Phi$ is a standard $(r,s)$-tensor field
on $\f M$ and $\nabla\Phi$ coincides with the standard covariant
differential with respect to the classical connection $\Gam$.
\END\eRm

\smallskip

\bCr\label{Cr4.4}
We have
\begin{align*}
\nab R[{K}] & =
        R[{K}]\col ji{\lam\mu;\nu}\, \E b^j\ten \E b_i
                \ten d^\lam\wed d^\mu\ten d^\nu
\\
& = \big( \der_{\nu}R[K]\col ji{\lam\mu} - K\col pi\nu\,
        R[K]\col jp{\lam\mu} + K\col jp\nu \,
        R[K]\col pi{\lam\mu}
\\
&\quad  + \Gam\col \nu\rho\lam \, R[K]\col ji{\rho\mu}
        + \Gam\col \nu\rho\mu\, R[K] \col ji{\lam\rho}\big)\,
        \E b^j\ten \E b_i \ten d^\lam\wed d^\mu\ten d^\nu\,.
        \END
\end{align*}
\eCr

The generalized Bianchi identity can be expressed by covariant
differentials as follows.

\bTh\label{Th4.5}
(The generalized Bianchi identity) We have
\beq
R[{K}]\col {j}i{\lam\mu;\nu}+
        R[{K}]\col {j}i{\mu\nu;\lam}+
        R[{K}]\col {j}i{\nu\lam;\mu}=0\,.
\eeq
\eTh

\bPf
This can be proved easily in coordinates by using
Corollary \ref{Cr4.4}. \QED
\ePf

\bTh\label{Th4.6}
Let $\Phi\in C^\infty(\f E^{p,r}_{q,s})$.
Then we have
\beq
\Alt \nab^2 \Phi =-\frac12\, R[\Gam^p_q \ten{K}^r_s]\com \Phi
        \in C^\infty(\f E^{p,r}_{q,s}\ten\Lam^2T^*\f M)\,,
\eeq
where $\Alt$ is the antisymmetrization.
\eTh

\bPf
This can be proved in coordinates by using
Proposition \ref{Pr3.5} and Proposition \ref{Pr4.3}. \QED
\ePf

\bEx\label{Ex4.7}
Let $\Phi\in C^\infty(\f E)$, $\Phi=\phi^i\E b_i$. Then
\beq
\Alt \nab^2 \Phi = - \frac12  R[{K}]\com \Phi
        :\f M\to \f E\ten \Lam^2T^*\f M\,,
\eeq
i.e. in coordinates
\beq
\Alt \nab^2 \Phi = -\frac12\, R[{K}]\col ji{\lam\mu}\, \phi^j \,
        \E b_i\ten d^\lam\wed d^\mu \,.\quad \END
\eeq
\eEx

\bEx\label{Ex4.8}
We have
\beq
\Alt \nab^2 R[{K}] :\f M\to\f E^*\ten\f E \ten
        \Lam^2T^*\f M\ten\Lam^2T^*\f M
        \,,
\eeq
expressed in coordinates by
\begin{align}\nonumber
\Alt \nab^2 R[{K}] & = -\frac12 \, \big(R[K]_{p}{}^i{}_{\nu_1\nu_2}\,
        R[{K}]_j{}^p{}_{\lam\mu} -R[K]_{j}{}^p{}_{\nu_1\nu_2}\,
        R[{K}]_p{}^i{}_{\lam\mu}
\\
&\quad -
        R[\Gam]_{\lam}{}^\ome{}_{\nu_1\nu_2} \,
        R[{K}]_j{}^i{}_{\ome\mu}
        - R[\Gam]_{\mu}{}^\ome{}_{\nu_1\nu_2} \,
        R[{K}]_j{}^i{}_{\lam\ome}\big) \,\nonumber
\\
& \qquad
\E b^j\ten \E b_i\ten d^\lam\wed d^\mu\ten d^{\nu_1}\wed d^{\nu_2}
        \,.\,\,\, \END\nonumber
\end{align}
\eEx



\end{document}